\documentclass[a4paper,11pt]{article}
\usepackage[colorinlistoftodos,prependcaption,textsize=tiny]{todonotes}

\usepackage{graphicx}
\usepackage{amssymb,amsmath,amsthm,times,mathrsfs,fullpage}
\usepackage{tikz,wrapfig,multirow}
\usetikzlibrary{patterns,arrows,decorations.pathreplacing}
\usepackage[pdftex]{hyperref}
\usepackage{epstopdf}
\hypersetup{plainpages=True, pdfstartview=FitH, bookmarksopen=true, colorlinks=true,linkcolor=blue,citecolor=blue}

\def\pf{\noindent \emph{Proof.}\ }
\def\qed{{\quad\rule{1mm}{3mm}\,}}

\begin{document}

\pgfdeclarelayer{background}
\pgfdeclarelayer{foreground}
\pgfsetlayers{background,main,foreground}

\newtheorem{thm}{Theorem}
\newtheorem{cor}{Corollary}
\newtheorem{lmm}{Lemma}
\newtheorem{conj}{Conjecture}
\newtheorem{pro}{Proposition}
\newtheorem{Def}{Definition}
\theoremstyle{remark}\newtheorem{Rem}{Remark}

\title{Asymptotic Enumeration and Distributional Properties of \\ Galled Networks}
\author{Michael Fuchs\thanks{MF was partially supported by grant MOST-109-2115-M-004-003-MY2.}\\
    Department of Mathematical Sciences\\
    National Chengchi University\\
    Taipei 116\\
    Taiwan \and
    Guan-Ru Yu\\
    Department of Mathematics\\
    National Kaohsiung Normal University\\
    Kaohsiung 824\\
    Taiwan \and
    Louxin Zhang\thanks{LZ was financially supported by Singapore MOE Tier-1 Research Fund R-146-000-318-114.}\\
    Department of Mathematics\\
    National University of Singapore\\
    Singapore 119076\\
    Republic of Singapore}

\maketitle

\begin{abstract}
We show a first-order asymptotics result for the number of galled networks with $n$ leaves. This is the first class of phylogenetic networks of {\it large} size for which an asymptotic counting result of such strength can be obtained. In addition, we also find the limiting distribution of the number of reticulation nodes of a galled networks with $n$ leaves chosen uniformly at random. These results are obtained by performing an asymptotic analysis of a recent approach of Gunawan, Rathin, and Zhang (2020) which was devised for the purpose of (exactly) counting galled networks. Moreover, an old result of Bender and Richmond (1984) plays a crucial role in our proofs, too.
\end{abstract}

\section{Introduction and Results}\label{intro}

Over the last few decades, {\it phylogenetic networks} have become a fundamental tool in evolutionary biology. Their now wide-spread usage makes it necessary to understand their basic combinatorial properties such as counting them or understanding the distribution of shape parameters when they are picked uniformly at random. Several recent papers have been dedicated to such studies; see \cite{BiLaSt, BoGaMa, CaZh, DiSeWe, FuGiMa, FuGiMa2, GuRaZh, Zh1}. The goal of this paper is to prove asymptotic counting results and investigate the stochastic behavior of shape parameters of {\it galled networks} which will be defined below.

First, a (binary, rooted) phylogenetic network with $n$ leaves is defined as connected, rooted  directed acyclic graph (DAG) whose nodes can be classified into four categories:
\begin{itemize}
\item[(a)] a {\it root} of indegree $0$ and outdegree $1$;
\item[(b)] {\it leaves} which are nodes of indegree $1$ and outdegree $0$ and which are bijectively labeled with $1,\ldots, n$;
\item[(c)] {\it tree nodes} which are nodes of indegree $1$ and outdegree $2$;
\item[(d)] {\it reticulation nodes} which are nodes of indegree $2$ and outdegree $1$.
\end{itemize}

Note that there are infinite phylogenetic networks with $n$ leaves for $n\geq 2$. However,  many interesting subclasses of phylogenetic networks contain a finite number of networks, e.g., normal networks, tree-child networks, galled networks, reticulation-visible networks, etc.; for definitions see, e.g., \cite{Book_Zhang2} and the references therein.

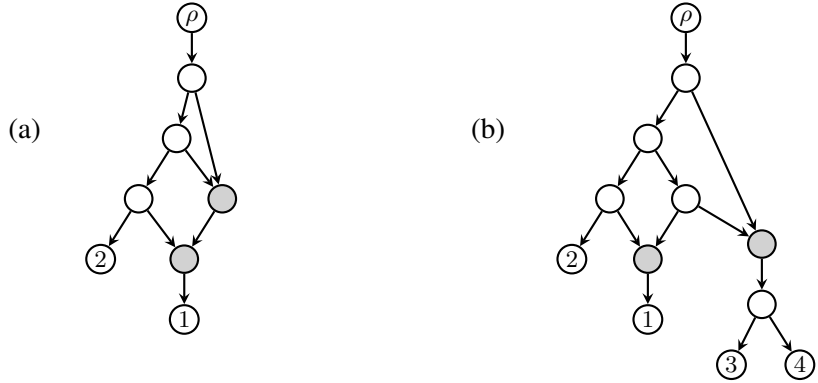
\begin{figure}[h]
\begin{center}
\begin{tikzpicture}
\draw (-2.2cm,-1.5cm) node (0) {(a)};
\draw (0cm,0cm) node[inner sep=1.2pt,line width=0.8pt,draw,circle] (1) {{\footnotesize $\rho$}};
\draw (0cm,-0.8cm) node[minimum size=7.8pt,line width=0.8pt,draw,circle] (2) {};
\draw (-0.2cm,-1.6cm) node[minimum size=7.8pt,line width=0.8pt,draw,circle] (3) {};
\draw (-0.7cm,-2.4cm) node[minimum size=7.8pt,line width=0.8pt,draw,circle] (4) {};
\draw (0.4cm,-2.4cm) node[minimum size=7.8pt,line width=0.8pt,draw,circle,fill=gray!35] (5) {};
\draw (-1.2cm,-3.2cm) node[inner sep=1.2pt,line width=0.8pt,draw,circle] (7) {{\footnotesize $2$}};
\draw (-0.1cm,-3.2cm) node[minimum size=7.8pt,line width=0.8pt,draw,circle,fill=gray!35] (8) {};
\draw (-0.1cm,-4cm) node[inner sep=1.2pt,line width=0.8pt,draw,circle] (9) {{\footnotesize $1$}};

\draw[-stealth,line width=0.8pt] (1) -- (2);
\draw[-stealth,line width=0.8pt] (2) -- (3);
\draw[-stealth,line width=0.8pt] (2) -- (5);
\draw[-stealth,line width=0.8pt] (3) -- (4);
\draw[-stealth,line width=0.8pt] (3) -- (5);
\draw[-stealth,line width=0.8pt] (4) -- (8);
\draw[-stealth,line width=0.8pt] (5) -- (8);
\draw[-stealth,line width=0.8pt] (4) -- (7);
\draw[-stealth,line width=0.8pt] (8) -- (9);

\draw (3.9cm,-1.5cm) node (0) {(b)};
\draw (6.5cm,0cm) node[inner sep=1.2pt,line width=0.8pt,draw,circle] (1) {{\footnotesize $\rho$}};
\draw (6.5cm,-0.8cm) node[minimum size=7.8pt,line width=0.8pt,draw,circle] (2) {};
\draw (6cm,-1.6cm) node[minimum size=7.8pt,line width=0.8pt,draw,circle] (3) {};
\draw (7.5cm,-3cm) node[minimum size=7.8pt,line width=0.8pt,draw,circle,fill=gray!35] (6) {};
\draw (5.5cm,-2.4cm) node[minimum size=7.8pt,line width=0.8pt,draw,circle] (4) {};
\draw (6.5cm,-2.4cm) node[minimum size=7.8pt,line width=0.8pt,draw,circle] (5) {};
\draw (5cm,-3.2cm) node[inner sep=1.2pt,line width=0.8pt,draw,circle] (7) {{\footnotesize $2$}};
\draw (6cm,-3.2cm) node[minimum size=7.8pt,line width=0.8pt,draw,circle,fill=gray!35] (8) {};
\draw (6cm,-4cm) node[inner sep=1.2pt,line width=0.8pt,draw,circle] (9) {{\footnotesize $1$}};
\draw (7.5cm,-3.8cm) node[minimum size=7.8pt,line width=0.8pt,draw,circle] (10) {};
\draw (7.1cm,-4.6cm) node[inner sep=1.2pt,line width=0.8pt,draw,circle] (11) {{\footnotesize $3$}};
\draw (8cm,-4.6cm) node[inner sep=1.2pt,line width=0.8pt,draw,circle] (12) {{\footnotesize $4$}};

\draw[-stealth,line width=0.8pt] (1) -- (2);
\draw[-stealth,line width=0.8pt] (2) -- (3);
\draw[-stealth,line width=0.8pt] (2) -- (6);
\draw[-stealth,line width=0.8pt] (3) -- (4);
\draw[-stealth,line width=0.8pt] (3) -- (5);
\draw[-stealth,line width=0.8pt] (4) -- (8);
\draw[-stealth,line width=0.8pt] (5) -- (8);
\draw[-stealth,line width=0.8pt] (5) -- (6);
\draw[-stealth,line width=0.8pt] (4) -- (7);
\draw[-stealth,line width=0.8pt] (8) -- (9);
\draw[-stealth,line width=0.8pt] (6) -- (10);
\draw[-stealth,line width=0.8pt] (10) -- (11);
\draw[-stealth,line width=0.8pt] (10) -- (12);
\end{tikzpicture}
\end{center}
\caption{(a) A non-galled network, in which the reticulation node at the bottom is not in any tree cycle. (b) A galled network.}\label{galled-net}
\end{figure}

Here, we will investigate the counting-related issues of {\it galled networks}. In a phylogenetic network, a {\it tree cycle} is a union of two edge-disjoint paths that are from a tree node to  a reticulation node with all other nodes being tree nodes. Galled networks are defined as phylogenetic networks with every reticulation node contained in a tree cycle; see Figure~\ref{galled-net} for examples.

\vspace*{0.35cm}

The above listed subclasses of phylogenetic networks are all proved or expected to be significantly ``larger" than the class of binary phylogenetic trees. For instance, for normal and tree-child networks, their sizes are  proved  to  grow (up to smaller-order terms in the main asymptotics) like $n^{2n}$ (see \cite{DiSeWe}). The numbers of the other subclasses above are expected to grow at a similar large speed. This is in contrast to, e.g., level $1$ and level $2$ networks or normal and tree-child networks with a fixed number of reticulation nodes which all grow like $n^{n}$ (see \cite{BoGaMa, FuGiMa, FuGiMa2}). Note that this is the same speed of growths as exhibited by binary phylogenetic trees with $n$ leaves; see, e.g., \cite{DiSeWe}. For all these latter ``small" subclasses, even the precise first-order asymptotics for their numbers is known; see also \cite{BoGaMa, FuGiMa, FuGiMa2, DiSeWe}.

On the other hand, for ``large" subclasses, the first-order asymptotics is still unknown for all major classes even though we came quite close in establishing a first-order asymptotic result for the number of tree-child networks in \cite{FuYuZh}. In order to recall our result, denote by $\mathrm{TC}_n$ the number of tree-child networks with $n$ leaves. Then, we proved in \cite{FuYuZh} that, as $n\rightarrow\infty$,
\[
\mathrm{TC}_n=\Theta\left(n^{-2/3}e^{a_1(3n)^{1/3}}\left(\frac{12}{e^2}\right)^nn^{2n}\right),
\]
where $a_1$ is the largest root of the Airy function of first kind.

The surprise here was the presence of the {\it stretched exponential} $e^{a_1(3n)^{1/3}}$ in the asymptotics. In the conclusion of \cite{FuYuZh}, we asked whether such a stretched exponential is also present in the asymptotics of other ``large" subclasses of phylogenetic networks.

In this paper, we show that this guess is wrong for the number of galled networks with $n$ leaves which we are going to denote by $\mathrm{GN}_n$. In fact, this class is the first ``large" class for which we are able to find a precise first-order asymptotic result.

\begin{thm}\label{main-thm-1}
For the number $\mathrm{GN}_n$ of galled networks with $n$ leaves, as $n\rightarrow\infty$,
\[
\mathrm{GN}_n\sim\frac{\sqrt{2e\sqrt[4]{e}}}{4}n^{-1}\left(\frac{8}{e^2}\right)^nn^{2n}.
\]
\end{thm}

As an important intermediate step in the proof of this result, we will also derive the first-order asymptotics of the number of {\it one-component galled networks} which are galled networks where the child of every reticulation node is a leaf. They are used as a building block in the construction of (general) galled networks; see \cite{GuRaZh} and the next section. We will denote the number of one-component galled networks with $n$ leaves by $\mathrm{OGN}_n$ throughout this work.

\begin{pro}\label{main-pro-1}
For the number $\mathrm{OGN}_n$ of one-component galled networks with $n$ leaves, as $n\rightarrow\infty$,
\[
\mathrm{OGN}_n\sim\frac{\sqrt{2\sqrt{e}}}{4}n^{-1}\left(\frac{8}{e^2}\right)^nn^{2n}.
\]
In particular, the fraction of one-component galled networks with $n$ leaves amongst (general) galled networks with $n$ leaves tends to $e^{-3/8}$ as $n$ tends to infinity.
\end{pro}

\begin{Rem}
\begin{itemize}
\item[(a)] The first-order asymptotics of the number of one-component tree-child networks (defined in a similar way as one-component galled networks and denoted by $\mathrm{OTC}_n$) was found in \cite{FuYuZh}, which showed that $\mathrm{OTC}_n$ is smaller than $\mathrm{TC}_n$ by an exponential order. Thus, in contrast to galled networks, the fraction of one-component tree-child networks with $n$ leaves amongst (general) tree-child networks with $n$ leaves tends to $0$ exponentially fast as $n$ tends to infinity.
\item[(b)] It was proved in \cite{CaZh} that the classes of one-component galled, reticulation-visible, tree-based, and phylogenetic networks all collapse into one. Thus, our above result gives the first-order asymptotics for all these classes of one-component networks, too.
\end{itemize}
\end{Rem}

Our approach for proving Theorem~\ref{main-thm-1} (and Proposition~\ref{main-pro-1}) will also allow us, for the first time, to give a detailed study of shape parameters of random networks, where the random model is the uniform model, i.e., networks are picked with identical probabilities.

We first need some notations. We call a reticulation node of a galled network {\it inner} if its child is not a leaf. So, for instance, one-component galled networks have no inner reticulation nodes, whereas the galled network in  Figure~\ref{galled-net}-(b) has exactly one inner reticulation node.

Next, denote by $X_n$ the number of inner reticulation nodes of a random galled network with $n$ leaves and by $Y_n$ the total number of reticulation nodes. Then, we have the following limit distribution result.

\begin{thm}\label{main-thm-2}
The random vector $(X_n,n-Y_n)$ weakly tends to a discrete limit distribution $(X,Y)$, i.e., as $n\rightarrow\infty$,
\[
(X_n,n-Y_n)\stackrel{d}{\longrightarrow}(X,Y),
\]
where $\stackrel{d}{\longrightarrow}$ denotes convergence in distribution. Moreover, the limit law of $(X,Y)$ is given by
\[
P(X=j,Y=k)=\frac{e^{-7/8}}{16^jj!}[z^{j-k}]e^{1/(2z)}(1+2z+3z^2)^j,\qquad (j\geq 0, k\geq -j),
\]
where $[z^{n}]f(z)$ denotes the $n$-th coefficient in the power series expansion of $f(z)$ centered at $0$.
\end{thm}

This result has two consequences.

\begin{cor}\label{first-cor}
\begin{itemize}
\item[(i)] The number of reticulation nodes $Z_n$ of a one-component galled network with $n$ leaves picked uniformly at random satisfies the following limit distribution result:
\[
n-Z_n\stackrel{d}{\longrightarrow}\mathrm{Poi}(1/2),\qquad (n\rightarrow\infty),
\]
where $\mathrm{Poi}(\lambda)$ denotes the Poisson distribution with parameter $\lambda$.
\item[(ii)] The limit distribution $X$ of the number of inner reticulation nodes $X_n$ of a galled network with $n$ leaves picked uniformly at random is a Poisson distribution with parameter $3/8$.
\end{itemize}
\end{cor}

\begin{cor}\label{second-cor}
The mean and the variance of the number of reticulation nodes $Y_n$ of a galled network with $n$ leaves picked uniformly at random satisfies, as $n\rightarrow\infty$,
\[
{\mathbb E}(Y_n)=n-\frac{3}{8}+o(1)\qquad\text{and}\qquad\mathrm{Var}(Y_n)=\frac{3}{4}+o(1).
\]
\end{cor}

The proofs of all the above results will rest on a recent approach developed by Guanawan et al. \cite{GuRaZh} to count exactly galled networks.
Since the approach is crucial in our asymptomic analyses, we will recap it in the next section. Moreover, we will also give tables for the counts when the number $n$ of leaves is small with some of the entries correcting those of the corresponding tables from \cite{GuRaZh}. In Section~\ref{prep}, we will show inequalities and discuss monotonicity properties which will turn out to be important in the proof of Proposition~\ref{main-pro-1}. In Section~\ref{enum-1-c}, we will prove Proposition~\ref{main-pro-1}, i.e., derive the first-order asymptotics of the number of one-component galled networks. Roughly speaking, this result will be deduced from a recurrence given in \cite{GuRaZh} and the proof will eventually rest on the Laplace method; see Appendix B.6 in \cite{FlSe} or Chapter~9 in \cite{GrKnPa}. In Section~\ref{enum-galled}, we will prove Theorem~\ref{main-thm-1}, i.e., derive the first-order asymptotics of the number of (general) galled networks, by applying the following strategy. First, we will use a result from \cite{BeRi} to obtain an asymptotic upper bound. Then, with the help of this upper bound, we will be able to identify the galled networks which asymptotically dominate. Finally, we will derive the first-order asymptotics of the number of these networks giving a matching lower bound. The insights from Section~\ref{enum-galled} will also be crucial for the proof of Theorem~\ref{main-thm-2} which will be presented, together with the proofs of its corollaries,  in Section~\ref{num-of-ret}.
Finally, in Section~\ref{dup-trees}, we will explain that our approach can also be used to count dup-trees which are leaf-labeled trees where every label from the label set $\{1,\ldots,n\}$ can be used at most twice. We will finish the paper with some concluding remarks in Section~\ref{con}.

\section{The Approach of Gunawan, Rathin, and Zhang for Counting Galled Networks}

The purpose of this section is to explain the approach for (exact) enumeration of galled networks that appeared in \cite{GuRaZh}.

First, consider one-component galled networks. If we denote by $\mathrm{OGN}_{n,k}$ the number of one-component galled networks with $n$ leaves and $k$ reticulation nodes, we have:
\begin{equation}
\label{rel-GNn-GNnk}
\mathrm{OGN}_n=\sum_{k=0}^{n}\mathrm{OGN}_{n,k}.
\end{equation}
Also, if $t$ denotes the number of tree nodes, then an easy counting argument shows that
\begin{equation}\label{rel-nkt}
n+k=t+1;
\end{equation}
see \cite{DiSeWe}. (This in fact holds for any phylogenetic network.) Consequently, the number of one-component galled networks with $n$ leaves is finite.

The (exact) counting problem for one-component galled networks with $n$ leaves was solved in \cite{GuRaZh}. More precisely, the authors of \cite{GuRaZh} proved the following result.

\begin{pro}[Gunawan et al. \cite{GuRaZh}]\label{gun-one-comp}
The number $\mathrm{OGN}_{n,k}$ of one-component galled networks with $n$ leaves and $k$ reticulation nodes is given by
\[
\mathrm{OGN}_{n,k}=\binom{n}{k}N_{n+1}^{(k)},
\]
where $N_{n}^{(k)}$ denotes the number of one-component galled networks with $n-1$ leaves and with the children of the reticulation nodes being labeled with $1,\dots,k$. Moreover, this number is recursively given by
\begin{align}
N_{n}^{(k)}=(n+k-3)N_{n}^{(k-1)}&+(k-1)N_{n}^{(k-2)}\nonumber\\
&+\frac{1}{2}\sum_{1\leq d\leq k-1}\binom{k-1}{d}(2d-1)!!\left(N_{n-d}^{(k-1-d)}-N_{n-d+1}^{(k-1-d)}\right)\label{rec-Nki}
\end{align}
for $2\leq k\leq n-1$ with initial values $N_{n}^{(0)}=(2n-5)!!$ and $N_n^{(1)}=(n-2)(2n-5)!!$.
\end{pro}

Using the recurrence in Proposition~\ref{gun-one-comp}, $N_{n}^{(k)}$ for small values of $n$ and $k$ can be computed; see Table~\ref{n-n-k} in the appendix which coincides with Table~1 in \cite{GuRaZh} with the only difference that some of the entries are corrected (namely, the value for $n=10$ and $k=9$ and the values for $n=11$ and $7\leq k\leq 10$).

\begin{figure}[h]
\begin{center}
\includegraphics[scale=0.35]{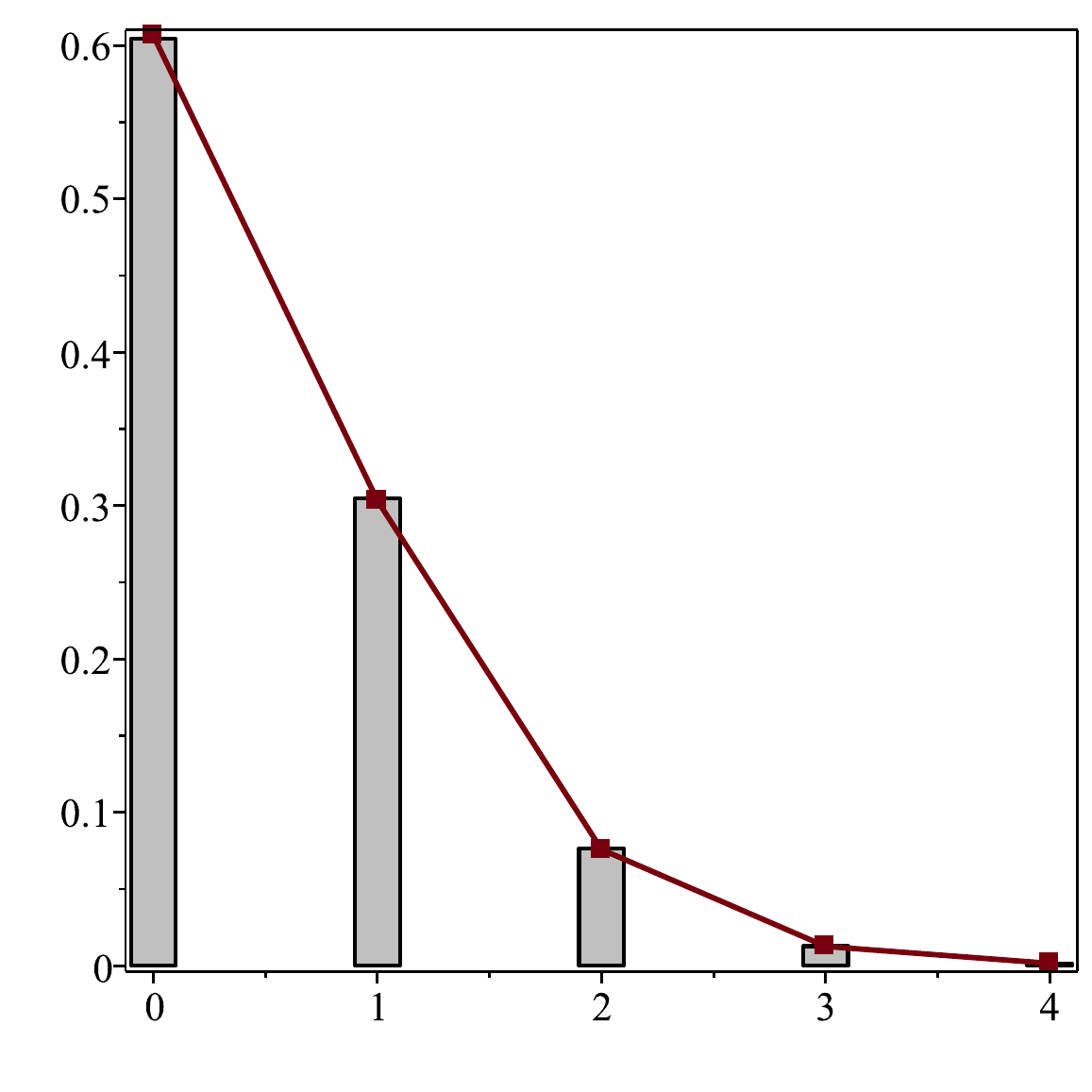}
\end{center}
\vspace*{-0.5cm}
\caption{A plot of $P(n-Z_{n}=k)$ with $n=100$ and $k=0, 1, 2, 3, 4$ (gray bars) and the claimed limit distribution (${\rm Poi}(1/2)$; filled squares) from Corollary~\ref{first-cor}.}\label{histo}
\end{figure}

Moreover, the values of $\mathrm{OGN}_{n,k}$ and $\mathrm{OGN}_{n}$ can be computed as well and thus also the distribution of $n-Z_n$ since
\[
{\mathbb P}(n-Z_n=k)=\frac{\mathrm{OGN}_{n,n-k}}{\mathrm{OGN}_{n}};
\]
see Figure~\ref{histo} for a plot of the histogram of $n-Z_n$ for $n=100$ from which the claimed limit law of Corollary~\ref{first-cor}-(i) is visible.

\begin{figure}[tp]
\begin{center}
\begin{tikzpicture}
\draw (0cm,0cm) node[inner sep=1.2pt,line width=0.8pt,draw,circle] (1) {{\footnotesize $\rho$}};
\draw (0cm,-0.8cm) node[minimum size=7.8pt,line width=0.8pt,draw,circle] (2) {};
\draw (-.5cm,-1.6cm) node[minimum size=7.8pt,line width=0.8pt,draw,circle] (3) {};
\draw (1cm,-3cm) node[minimum size=7.8pt,line width=0.8pt,draw,circle,fill=gray!35] (6) {};
\draw (-1cm,-2.4cm) node[minimum size=7.8pt,line width=0.8pt,draw,circle] (4) {};
\draw (0cm,-2.4cm) node[minimum size=7.8pt,line width=0.8pt,draw,circle] (5) {};
\draw (-1.5cm,-3.2cm) node[inner sep=1.2pt,line width=0.8pt,draw,circle] (7) {{\footnotesize $2$}};
\draw (-0.5cm,-3.2cm) node[minimum size=7.8pt,line width=0.8pt,draw,circle,fill=gray!35] (8) {};
\draw (-0.5cm,-4cm) node[inner sep=1.2pt,line width=0.8pt,draw,circle] (9) {{\footnotesize $1$}};
\draw (1cm,-3.8cm) node[minimum size=7.8pt,line width=0.8pt,draw,circle] (10) {};
\draw (0.6cm,-4.6cm) node[inner sep=1.2pt,line width=0.8pt,draw,circle] (11) {{\footnotesize $3$}};
\draw (1.5cm,-4.6cm) node[inner sep=1.2pt,line width=0.8pt,draw,circle] (12) {{\footnotesize $4$}};

\draw[-stealth,line width=0.8pt] (1) -- (2);
\draw[-stealth,line width=0.8pt] (2) -- (3);
\draw[-stealth,line width=0.8pt] (2) -- (6);
\draw[-stealth,line width=0.8pt] (3) -- (4);
\draw[-stealth,line width=0.8pt] (3) -- (5);
\draw[-stealth,line width=0.8pt] (4) -- (8);
\draw[-stealth,line width=0.8pt] (5) -- (8);
\draw[-stealth,line width=0.8pt] (5) -- (6);
\draw[-stealth,line width=0.8pt] (4) -- (7);
\draw[-stealth,line width=0.8pt] (8) -- (9);
\draw[-stealth,line width=0.8pt] (6) -- (10);
\draw[-stealth,line width=0.8pt] (10) -- (11);
\draw[-stealth,line width=0.8pt] (10) -- (12);

\draw[stealth-stealth,line width=0.8pt] (3cm,-2.3cm) -- (5cm,-2.3cm);

\draw (7.5cm,-1.3cm) node[minimum size=7.8pt,line width=0.8pt,draw,circle] (1) {};
\draw (6.7cm,-2.3cm) node[inner sep=1.2pt,line width=0.8pt,draw,circle] (2) {{\footnotesize $2$}};
\draw (7.5cm,-2.3cm) node[inner sep=1.2pt,line width=0.8pt,draw,circle] (3) {{\footnotesize $1$}};
\draw (8.3cm,-2.3cm) node[minimum size=7.8pt,line width=0.8pt,draw,circle] (4) {};
\draw (7.9cm,-3.3cm) node[inner sep=1.2pt,line width=0.8pt,draw,circle] (5) {{\footnotesize $3$}};
\draw (8.7cm,-3.3cm) node[inner sep=1.2pt,line width=0.8pt,draw,circle] (6) {{\footnotesize $4$}};

\draw [line width=0.8pt] (1) -- (2);
\draw [line width=0.8pt] (1) -- (3);
\draw [line width=0.8pt] (1) -- (4);
\draw [line width=0.8pt] (4) -- (5);
\draw [line width=0.8pt] (4) -- (6);

\filldraw (7.42cm,-1.75cm) -- (7.58cm,-1.75cm) -- (7.5cm,-1.91cm);
\filldraw (7.92cm,-1.7cm) -- (7.8cm,-1.8cm) -- (7.97cm,-1.89cm);
\end{tikzpicture}
\end{center}
\caption{A galled network and its underlying phylogenetic tree where the triangle on an edge indicates that the head of the edge is a reticulation node in the original galled network.}\label{galled-net-decom}
\end{figure}
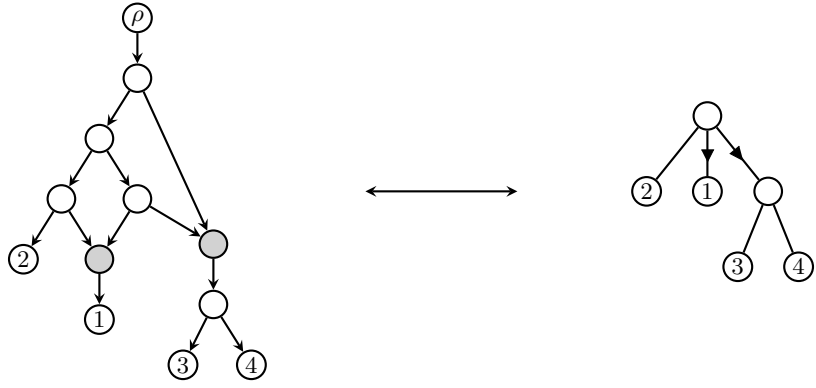

Next, we consider (general) galled networks. Here, the crucial idea of the approach used in \cite{GuRaZh} is the decomposition of phylogenetic networks into tree-node components \cite{Book_Zhang2}. More precisely, each  galled network can be compressed into a rooted phylogenetic  tree (with all non-leaf nodes having at least two children)  if  its tree-node components are replaced with single nodes;  see Figure~\ref{galled-net-decom} for the compression  of the galled network in Figure~\ref{galled-net}-(b). Note that each tree-node component together with all incident reticulation nodes form a one-component galled network.

This decomposition is also reversible, i.e., every galled network can be constructed by starting with a phylogenetic tree and then replacing all non-leaf nodes by one-component galled networks whose number of leaves is equal to the outdegree of the replaced node. In addition, after the replacement of a non-leaf node, its children which have been internal nodes must be below reticulation nodes and its children which have been leaves may or may not be below reticulation nodes; see again Figure~\ref{galled-net-decom} where arrows on edges to children indicate that they are below reticulation nodes when their parent is replaced.

As a result of the above procedure, the following formula for computing the number of galled networks with $n$ leaves was given in \cite{GuRaZh}.

\begin{thm}[Gunawan et al. \cite{GuRaZh}]\label{main-thm-GuRaZh}
The number $\mathrm{GN}_n$ of galled networks with $n$ leaves is given by
\begin{equation}\label{form-GNn}
\mathrm{GN}_n=\sum_{\mathcal{T}}\prod_{v\in\mathcal{I}(\mathcal{T})}\sum_{j=c_{\mathrm{nlf}}(v)}^{c(v)}\binom{c_{\mathrm{lf}}(v)}{j-c_{\mathrm{nlf}}(v)}N^{(j)}_{c(v)+1},
\end{equation}
where the first sum runs over all phylogenetic trees $\mathcal{T}$, $\mathcal{I}(\mathcal{T})$ denotes the set of internal nodes of ${\mathcal T}$, $c(v)$ denotes the number of children of $v$ and $c_{\mathrm{nlf}}(v)$ and $c_{\mathrm{lf}}(v)$ denote the number of children which are non-leaf nodes and leaves of ${\mathcal T}$, respectively.
\end{thm}

From this result, the number of galled networks with $n$ leaves for small values of $n$ can be computed; see Table~\ref{GN-n},  which corrects the last two counts in Table~2 in \cite{GuRaZh}.

\vspace*{0.1cm}
\begin{table}[h]
\begin{center}
\begin{tabular}{c|c}
$n$ & $\mathrm{GN}_n$ \\
\hline
1 & 1 \\
2 & 6 \\
3 & 240 \\
4 & 20,502 \\
5 & 2,868,990 \\
6 & 589,130,280 \\
7 & 167,357,180,970 \\
8 & 63,356,654,623,500 \\
9 & 31,092,212,800,634,580 \\
10 & 19,327,089,427,089,478,650 \\
\end{tabular}
\end{center}
\caption{The values of $\mathrm{GN}_n$ for $1\leq n\leq 10$.}\label{GN-n}
\end{table}

In fact, the above theorem also allows one to compute the number of galled networks with $n$ leaves and $k$ reticulation nodes which we are going to denote by $\mathrm{GN}_{n,k}$. We first make an easy observation.

\begin{lmm}\label{ret-max}
A galled network with $n$ leaves has at most $2n-2$ reticulation nodes and at most $n-2$ inner reticulation nodes.
\end{lmm}

\pf
Let $N$ be a galled network with $n$ leaves and $T$ be the underlying phylogenetic tree  of  the decomposition of $N$ into its tree-node components. $T$ has the same leaves as $N$. Assume $V_i(T)$ denotes the non-leaf nodes of $T$ and $d(v)$ denotes the number of the children of $v$ for each $v\in V_i(T)$. Since each non-root node has a unique parent and $T$ has $n$ leaves, $|V_i(T)|-1+n= \sum_{v\in V_i(T)}d(v)\geq 2|V_i(T)|$ and
thus  $V_i(T)\leq n-1$. Since every non-root node of $V_i(T)$ corresponds to an inner reticulation node of $N$,  $T$ has  at most $n-2$ inner reticulation nodes. Since each leaf of $T$  may or may not correspond to  a reticulation node of $N$, $N$ has at most $n+n-2$ reticulation nodes.

Moreover, if each non-leaf node of $T$ has exactly two children and if the parent of every leaf  is a reticulation node in $N$,  $N$  has $n-2$ inner reticulation nodes and $2n-2$ reticulation nodes. \qed

Using the above method of proof, we can also deduce the following formula for $\mathrm{GN}_{n,2n-2}$, i.e., for the number of galled networks with a maximal number of reticulation nodes.

\begin{lmm}
The number of galled networks with $n$ leaves and maximal number of reticulation nodes $2n-2$ is given by
\[
\mathrm{GN}_{n,2n-2}=\frac{(2n-2)!}{(n-1)!}\left(\frac{3}{2}\right)^{n-1}.
\]
\end{lmm}

\pf By the proof of Lemma~\ref{ret-max}, we see that galled networks with $n$ leaves and $2n-2$ reticulation nodes are constructed from binary phylogenetic trees with $n$ leaves by replacing all $n-1$ internal nodes by one-component galled networks with exactly two reticulation nodes. The number of binary phylogenetic trees with $n$ leaves is given by $\frac{(2n-2)!}{2^{n-1}(n-1)!}$ (see, e.g., Corollary 2.2.4 in \cite{SeSt}) and the number of one-component galled networks which replace internal nodes by $N_{3}^{(2)}=3$ (see \cite{GuRaZh}). Thus,
\[
\mathrm{GN}_{n,2n-2}=\frac{(2n-2)!}{2^{n-1}(n-1)!}\cdot 3^{n-1}
\]
which is the claimed result.\qed

With more work, $\mathrm{GN}_{n,2n-2-\ell}$ for small values of $\ell$ can be computed as well. Moreover, $\mathrm{GN}_{n,\ell}$ for small values of $\ell$ can also be found; see Proposition 21 in \cite{CaZh} and \cite{FuGiMa2} where such results were derived for normal and tree-child networks with $n$ leaves. However, we are less interested in such results because the most interesting range of $\mathrm{GN}_{n,k}$ will turn out to be $k$ close to $n$; compare with Theorem~\ref{main-thm-1}.

In fact, in order to compute the distribution of the random variable $(X_n,n-Y_n)$ from Theorem~\ref{main-thm-1} one needs the number of galled networks with $n$ leaves, $k$ reticulation nodes and $j$ inner reticulation nodes which we denote by $\mathrm{GN}_{n,k,j}$. Then,
\begin{equation}\label{joint-dist}
\mathbb{P}(X_n=j,n-Y_n=k-j)=\frac{\mathrm{GN}_{n,n+j-k,j}}{\mathrm{GN}_n},\qquad (0\leq j\leq n-2,0\leq k\leq n).
\end{equation}
These probabilities can be computed with the help of Theorem~\ref{main-thm-GuRaZh}, too; see Table~\ref{ll-ret} in the appendix for the values of the numerator of (\ref{joint-dist}) for $n=7$.

\section{Bounds and Monotonicity Properties}\label{prep}

In this section, we prove certain bounds for $N_{n}^{(k)}$ which then in turn imply bounds and montonicity of $\mathrm{OGN}_{n,k}$. These results will be needed in the next section for the asymptotic analysis of $\mathrm{OGN}_{n}$. As for $\mathrm{GN}_{n,k}$ and more generally $\mathrm{GN}_{n,k,j}$, we are not going to directly work with these sequences since we will be able to sidestep them in the proof of Theorem~\ref{main-thm-1} and Theorem~\ref{main-thm-2} thanks to a result in \cite{BeRi}. Thus, for these sequences, we will not need montonicity properties. Nevertheless, at the end of this section, we will give a brief discussion of which properties we expect for these sequences.

We start with a lower bound result for $N_{n}^{(k)}$.

\begin{lmm}\label{lb}
For $2\leq k\leq n-1$, we have $N_{n}^{(k)}\geq (n+k-3)N_{n}^{(k-1)}+\frac{k-1}{2}N_{n}^{(k-2)}$.
\end{lmm}
\pf Let ${\cal A}_{n-1,t}$ be the set of one-component galled networks with $n-1$ leaves and $t$ reticulation nodes whose children are labeled with $1, 2, \ldots, t$. Note that $N_{n}^{(t)}=|{\cal A}_{n-1, t}|$.

Let $G\in {\cal A}_{n-1,k-1}$. Then, $G$ contains $k-1$ reticulation nodes,  $(n-1)+(k-1)-1$ tree nodes and $n-1$ leaves (see (\ref{rel-nkt})). Since tree nodes and leaves are of indegree 1 and reticulation nodes are of indegree 2, $G$ contains $2(n-1)+3(k-1)-1=2n+3k-6$ edges (including the edge leaving the root).
For any edge $(x,y)$ such that $y$ is not equal to a leaf which is labeled with $i$ for $1\leq i\leq k$, inserting a node $u$ into $(x,y)$, inserting another node $v$ into $(p(k),k)$, where $p(k)$ denotes the parent of the leaf with label $k$, and adding the edge $(u,v)$ produces a one-component galled network $G'$ of ${\cal A}_{n-1, k}$. In this way, we can generate from $G$ $(2n+3k-6)-k=2(n+k-3)$ one-component galled networks of ${\cal A}_{n-1,k}$, some of which may be identical. Moreover, by different choices of $G$, we can generate
$2(n+k-3)N^{(k-1)}_n$ one-component galled networks of ${\cal A}_{n-1,k}$.

Next, for each $i$ such that $1\leq i\leq k-1$,  let ${\cal B}_{n-1,k-1,i}$ denote the set of one-component galled networks with $n-1$ leaves and $k-2$ reticulation nodes whose children have labels from $\{1,2,\ldots, k-1\}\setminus \{i\}$. Clearly, ${\cal B}_{n-1,k-1,i}$ contains $N^{(k-2)}_{n}$ networks. For each network $G \in {\cal B}_{n-1,k-1,i}$, the parents $p(i)$ and $p(k)$ of the leaves with labels $i$ and $k$ are tree nodes. Inserting $u$ and $v$ into $(p(i),i)$ to subdivide it into $(p(i),u)$, $(u,v)$ and $(v,i)$, inserting $w$ into $(p(k),k)$ and adding the edges $(u,t)$, $(t,v)$ and $(t,w)$ produces a network with $k$ reticulation nodes whose children are labeled with $1,2,\ldots k$. In this way, we can generate $\sum_{1\leq i\leq k-1}|{\cal B}_{n-1, k-1, i}|=(k-1)N^{(k-2)}_{n}$ networks of ${\cal A}_{n-1,k}$.

On the other hand, by removing either one of the edges entering the parent of the leaf with label $k$ from a network of ${\cal A}_{n-1, k}$ and contracting nodes of indegree 1 and outdegree 1 as well as double edges, we obtain:
\begin{itemize}
    \item a one-component galled network with $k-1$ reticulations whose children are leaves labeled with $1,2,\ldots,k-1$; or
    \item  a one-component galled network with $k-2$ reticulations whose children have labels from the set $\{1, 2, \cdots, k-1\}$.
\end{itemize}

Taken together, the two facts imply that $N_{n}^{(k)}\geq (n+k-3)N_{n}^{(k-1)}+\frac{k-1}{2}N_{n}^{(k-2)}$ which proves the claim.\qed

Next, we deduce from (\ref{rec-Nki}) an (almost matching) upper bound result for $N_n^{(k)}$.

\begin{lmm}\label{ub}
For $2\leq k\leq n-1$, we have $N_{n}^{(k)}\leq(n+k-3)N_{n}^{(k-1)}+\frac{k-1}{2}N_{n}^{(k-2)}+\frac{k-1}{2}N_{n-1}^{(k-2)}$.
\end{lmm}
\pf Note that for $0\leq k\leq n-2$, we have $N_{n}^{(k)}\geq N_{n-1}^{(k)}$. Thus, from (\ref{rec-Nki}), we obtain that
\begin{align*}
N_{n}^{(k)}&=(n+k-3)N_n^{(k-1)}+\frac{k-1}{2}N_{n}^{(k-2)}+\frac{k-1}{2}N_{n-1}^{(k-2)} \nonumber \\
   &\qquad\frac{1}{2}\sum_{2\leq d\leq k-1}\binom{k-1}{d}(2d-1)!!\left(N_{n-d}^{(k-1-d)}-N_{n-d+1}^{(k-1-d)}\right)  \nonumber \\
&\leq(n+k-3)N_n^{(k-1)}+\frac{k-1}{2}N_{n}^{(k-2)}+\frac{k-1}{2}N_{n-1}^{(k-2)}. \label{mod-rec-Nki}
\end{align*}
This proves the claimed result.\qed

The relation between $N_{n}^{(k)}$ and $N_{n-1}^{(k)}$ which was used in the above proof can actually be improved. (This improvement will be needed in the next section.)

\begin{lmm}\label{vert-bound}
For $0\leq k\leq n-2$, we have $N_{n}^{(k)}\geq(n+k-5/2)N_{n-1}^{(k)}$.
\end{lmm}
\pf Recall that a galled network with $n-2$ leaves and $k$ reticulation nodes has $2n+3k-5$ edges (see the proof of Lemma~\ref{lb}). Let $G$ be such a network. By choosing edges $(x,y)$ from $G$ such that $y$ is not a leaf with label $i$ for $1\leq i\leq k$, inserting a node $u$ and connecting $u$ to a new node with label $n-1$, we obtain $2n+2k-5$ one-component networks $G'$ with $n-1$ leaves and $k$ reticulation nodes. Conversely, note that each one-component network with $n-1$ leaves and $k$ reticulation nodes is obtained by the above procedure at most twice. From this, the claimed result follows.\qed

Now, we use the above results to prove corresponding bounds for $\mathrm{OGN}_{n,k}$.
First, from Lemma~\ref{lb}, we obtain the following.

\begin{cor}\label{prop_III}
Let $\mathrm{OGN}_{n,k}$ be the number of one-component galled networks with $n$ leaves and $k$ reticulation nodes. Then, we have the following facts:
\begin{itemize}
\item[(i)] For any $0\leq k\leq n-1$, we have $\mathrm{OGN}_{n,k+1}\geq\frac{(n-k)(n+k-1)}{k+1}\times\mathrm{OGN}_{n,k}$. Thus, $\mathrm{OGN}_{n,k}$ is increasing in $k$.
\item[(ii)] For any $0\leq k\leq n$, we have $\mathrm{OGN}_{n,k}\leq\frac{n!(n+k-2)!}{k!(n-k)!(2n-2)!}\times\mathrm{OGN}_{n,n}=\binom{n}{k}\frac{(n+k-2)!}{(2n-2)!}\times\mathrm{OGN}_{n,n}.$
\end{itemize}
\end{cor}
\pf (i) Recall that $\mathrm{OGN}_{n,k+1}=\binom{n}{k+1}N^{(k+1)}_{n+1}$. Moreover, from Lemma~\ref{lb}, we have
\begin{equation}\label{lbb}
N_{n}^{(k)}\geq (n+k-3)N_{n}^{(k-1)}
\end{equation}
for $2\leq k\leq n-1$ and this bound also holds for $k=1$. Now, we can estimate $\mathrm{OGN}_{n,k+1}$ as follows
\begin{align*}
 \mathrm{OGN}_{n,k+1}&=\binom{n}{k+1}N^{(k+1)}_{n+1}
  \geq \binom{n}{k+1}(n+k-1)N^{(k)}_{n+1}\\
   &\geq\frac{n-k}{k+1}\binom{n}{k} (n+k-1) N^{(k)}_{n+1}=\frac{(n-k)(n+k-1)}{k+1}\times\mathrm{OGN}_{n,k}.
 \end{align*}
For $n\geq 3$, we have $n+k-1>k+1$ and thus $\mathrm{OGN}_{n,k+1}>(n-k)\times\mathrm{OGN}_{n,k}$, i.e., $\mathrm{OGN}_{n, k}$ increases with $k$ for $n\geq 3$. Finally, that this property also holds for $n=2$ can easily be verified directly.

(ii) It can be derived from (i) by induction. \qed

Secondly, Lemma~\ref{vert-bound} implies the following relation, which can be proved in the same way as  Corollary~\ref{prop_III}.
\begin{cor}
For $0\leq k\leq n-1$, we have $\mathrm{OGN}_{n,k}\geq\frac{n(2n+2k-3)}{2(n-k)}\times\mathrm{OGN}_{n-1,k}.$
\end{cor}

Now, we will briefly discuss what we expect for the monotonicity behavior of $\mathrm{GN}_{n,k}$ and $\mathrm{GN}_{n,k,j}$; the claims below will not be needed in the sequel and proofs might appear elsewhere.

First consider $\mathrm{GN}_{n,k,j}$ where $j\leq k\leq n+j$. Note that $j=0$ corresponds to $\mathrm{OGN}_{n,k}$. We expect that $\mathrm{GN}_{n,k,j}$, when $j$ with $0\leq j\leq n-2$ is fixed, is increasing in the range $j\leq k\leq n$; compare with Table~\ref{ll-ret} for $n=7$. Since
\[
\mathrm{GN}_{n,k}=\sum_{j=0}^{\min\{n-2,k\}}\mathrm{GN}_{n,k,j},
\]
this would then imply that $\mathrm{GN}_{n,k}$ is increasing for $0\leq k\leq n$. On the other hand, the sequence $\mathrm{GN}_{n,k,j}$, again when $j$ is fixed, is in general not decreasing for $n<k\leq n+j$. In fact, the sequence seems to continue to increase for a few terms before it starts to decrease with the maximum more and more pushed to the right as $j$ gets large. However, since larger values of $j$ contribute less to $\mathrm{GN}_{n,k}$, we nevertheless expect that $\mathrm{GN}_{n,k}$ is decreasing for $n\leq k\leq 2n-2$. Overall, we have the following conjecture.

\begin{conj}
The sequence $\mathrm{GN}_{n,k}$ is increasing for $0\leq k\leq n$ and decreasing for $n\leq k\leq 2n-2$.
\end{conj}

This conjecture, if true, implies that the limit distribution $Y$ of the number of reticulation nodes from Theorem~\ref{main-thm-2} has a unique maximum at $0$ and the left and right tail are both non-increasing. (One might be able to deduce this directly from the expression of the limit distribution from Theorem~\ref{main-thm-2}.)

\section{Asymptotic Analysis of $\mathbf{OGN_{n}}$}\label{enum-1-c}

This section contains the proof of Proposition~\ref{main-pro-1}.

We start by deducing the following (asymptotic) simplification of the recurrence (\ref{rec-Nki}) for $N_{n}^{(k)}$ from Lemma~\ref{lb} and Lemma~\ref{ub}.

\begin{lmm}\label{approx-rec}
For $1\leq k\leq n-1$,
\begin{equation}
\label{key_formula}
N_{n}^{(k)}=\left(n+k-3+\frac{k}{2(n+k)}+{\mathcal O}\left(\frac{1}{n}\right)\right)N_{n}^{(k-1)},
\end{equation}
where the ${\mathcal O}$-estimate holds uniformly in $k$.
\end{lmm}
\pf  First note that from Lemma~\ref{lb} and Lemma~\ref{ub}, we have
\begin{equation}\label{first-est}
N_{n}^{(k)}=(n+k-3)N_{n}^{(k-1)}+\frac{k-1}{2}N_{n}^{(k-2)}+{\mathcal O}\left(kN_{n-1}^{(k-2)}\right),
\end{equation}
where the ${\mathcal O}$-estimate holds uniformly in $2\leq k\leq n-1$.

Next, again from Lemma~\ref{ub},
\[
N_n^{(k)}\leq (n+k-3)N_{n}^{(k-1)}+(k-1)N_{n}^{(k-2)}
\]
for $2\leq k\leq n-1$, since $N_{n}^{(k)}\geq N_{n-1}^{(k)}$. Consequently, from (\ref{lbb}),
\[
N_n^{(k)}\leq \left(n+k-3+\frac{k-1}{n+k-4}\right)N_{n}^{(k-1)}
\]
for $2\leq k\leq n-1$ and this also holds for $k=1$.

We use this now to estimate the second term on the right hand side of (\ref{first-est}):
\[
\frac{k-1}{2(n+k-4)+{\displaystyle\frac{2(k-2)}{n+k-5}}}N_{n}^{(k-1)}\leq\frac{k-1}{2}N_{n}^{(k-2)}\leq\frac{k-1}{2(n+k-4)}N_{n}^{(k-1)},
\]
where we again used (\ref{lbb}) for the upper bound. Note that
\[
\frac{k-1}{2(n+k-4)+{\displaystyle\frac{2(k-2)}{n+k-5}}}=\frac{k-1}{2(n+k-4)}\left(1+{\mathcal O}\left(\frac{1}{n}\right)\right)=\frac{k}{2(n+k)}+{\mathcal O}\left(\frac{1}{n}\right),
\]
where all ${\mathcal O}$-estimates hold uniformly in $2\leq k\leq n-1$. Consequently,
\begin{equation}\label{eqn_x6}
\frac{k-1}{2}N_{n}^{(k-2)}=\left(\frac{k-1}{2(n+k-4)}+{\mathcal O}\left(\frac{1}{n}\right)\right)N_{n}^{(k-1)}
\end{equation}
where the ${\mathcal O}$-estimate holds uniformly in $2\leq k\leq n-1$.

In addition, by Lemma~\ref{vert-bound} and again (\ref{lbb}),
\begin{equation}\label{eqn_x7}
kN_{n-1}^{(k-2)}\leq\frac{k}{2(n+k-9/2)} N_{n}^{(k-2)}\leq
\frac{k}{2(n+k-9/2)(n+k-4)} N_{n}^{(k-1)}={\mathcal O}\left(\frac{N_{n}^{(k-1)}}{n}\right),
\end{equation}
where again the ${\mathcal O}$-estimate holds uniformly in $2\leq k\leq n-1$.

Plugging (\ref{eqn_x6}) and (\ref{eqn_x7}) into (\ref{first-est}) implies (\ref{key_formula}). (Note that the case $k=1$ is trivial.) \qed

Now, we are ready to prove Proposition~\ref{main-pro-1} using the Laplace method. Recall that from (\ref{rel-GNn-GNnk}) and Proposition~\ref{gun-one-comp}, we have
\begin{equation}\label{rel-GNk-Nki}
\mathrm{OGN}_{n}=\sum_{k=0}^{n}\binom{n}{k}N_{n+1}^{(k)}.
\end{equation}

\vspace*{0.3cm}\noindent{\it Proof of Proposition~\ref{main-pro-1}.}
By Corollary~\ref{prop_III},  the terms in the sum (\ref{rel-GNk-Nki}) are increasing. Thus, the (asymptotic) main contribution to the sum is expected to come from large values of $k$.

Next, we consider the terms of (\ref{rel-GNk-Nki}) individually. By iterating the expression from Lemma~\ref{approx-rec} and using the following fact:
\[
n-1+\ell+\frac{\ell}{2(n+\ell)}+{\mathcal O}\left(\frac{1}{n}\right)=\left(n-1+\ell+\frac{\ell}{2(n+\ell)}\right)\left(1+{\mathcal O}\left(\frac{1}{n^2}\right)\right)
\]
uniformly for $0\leq\ell\leq n-1$, we obtain:
\begin{align}
N_{n+1}^{(k)}&= \left(\prod_{\ell=0}^{k-1}\left(n-1+\ell+\frac{\ell}{2(n+\ell)}+{\mathcal O}\left(\frac{1}{n}\right)\right)\right) N_{n+1}^{(0)}
 \nonumber \\
&=\left(\prod_{\ell=0}^{k-1}\left(n-1+\ell+\frac{\ell}{2(n+\ell)}\right)\right)\cdot\left(1+{\mathcal O}\left(\frac{1}{n^2}\right)\right)^{k-1} N_{n+1}^{(0)}, \label{ind-Nki}
\end{align}
where $0\leq k\leq n$ and $N_{n+1}^{(0)}=(2n-3)!!$.

The product term in the right-handed side of (\ref{ind-Nki}) can  further be simplified as:
\begin{align}
 \prod_{\ell=0}^{k-1}\Bigg(n -1&+\ell+\frac{\ell}{2(n-2+\ell)}\Bigg) \nonumber \\
&=\frac{(n-1)!{\displaystyle\Gamma\left(n+k-\frac{1}{4}-\frac{\sqrt{8n+1}}{4}\right)
\Gamma\left(n+k-\frac{1}{4}+\frac{\sqrt{8n+1}}{4}\right)}}{(n+k-1)!{\displaystyle\Gamma\left(n-\frac{1}{4}-\frac{\sqrt{8n+1}}{4}\right)
\Gamma\left(n-\frac{1}{4}+\frac{\sqrt{8n+1}}{4}\right)}}. \label{eqn15}
\end{align}

Plugging (\ref{eqn15}) into (\ref{ind-Nki}), multiplying by ${\displaystyle \binom{n}{k}}$ and replacing $k$ by $n-k$ gives
\[
\binom{n}{n-k}N_{n+1}^{(n-k)}=c_{n}^{(k)}\cdot\left(1+{\mathcal O}\left(\frac{1}{n^2}\right)\right)^{n-k-1}\cdot d_n,
\]
where, by an application of Stirling's formula,
\begin{align*}
c_n^{(k)}&=\frac{{\displaystyle\Gamma\left(2n-k-\frac{1}{4}-\frac{\sqrt{8n+1}}{4}\right)
\Gamma\left(2n-k-\frac{1}{4}+\frac{\sqrt{8n+1}}{4}\right)}}{k!(n-k)!(2n-k-1)!} \\
&=\frac{\sqrt[4]{e}}{k!2^{k+1}}n^{-3/2}\left(\frac{4}{e}\right)^nn^{n}\left(1+{\mathcal O}\left(\frac{k^2+1}{n}\right)\right)
\end{align*}
uniformly for $k=o(\sqrt{n})$,
and, similarly,
\begin{align*}
d_n&=\frac{n!(n-1)!(2n-3)!!}{\displaystyle\Gamma\left(n-\frac{1}{4}-\frac{\sqrt{8n+1}}{4}\right)
\Gamma\left(n-\frac{1}{4}+\frac{\sqrt{8n+1}}{4}\right)}.\\
 &= \frac{1}{\sqrt{2e}}n^{1/2}\left(\frac{2}{e}\right)^{n}n^{n}\left(1+{\mathcal O}\left(\frac{1}{n}\right)\right).
\end{align*}

Finally, since
\[
\left(1+{\mathcal O}\left(\frac{1}{n^2}\right)\right)^{n-k-1}
\leq \left(1+{\mathcal O}\left(\frac{1}{n^2}\right)\right)^{n} =
e^{{\mathcal O}(1/n)}=1+{\mathcal O}\left(\frac{1}{n}\right)
\]
for any $k$ such that $0\leq k\leq n$,
\begin{equation}\label{asymp-binom-Nk}
\binom{n}{n-k}N_{n+1}^{(n-k)}=\frac{\sqrt{2}}{k!2^{k+2}\sqrt[4]{e}}n^{-1}\left(\frac{8}{e^2}\right)^n n^{2n}\left(1+{\mathcal O}\left(\frac{k^2+1}{n}\right)\right)
\end{equation}
uniformly in $k=o(\sqrt{n})$.

Now, we break the sum in (\ref{rel-GNk-Nki}) into two parts:
\[
\mathrm{OGN}_n=\sum_{k<\sqrt[4]{n}}\binom{n}{n-k}N_{n+1}^{(n-k)}+\sum_{k\geq\sqrt[4]{n}}\binom{n}{n-k}N_{n+1}^{(n-k)}=:\Sigma_{1}+\Sigma_{2}.
\]
For the first sum by using (\ref{asymp-binom-Nk}),
\begin{align*}
\Sigma_1&=\left(\sum_{k<\sqrt[4]{n}}\frac{1}{k!2^{k}}\right)\cdot\frac{\sqrt{2}}{4\sqrt[4]{e}}n^{-1}\left(\frac{8}{e^2}\right)^n n^{2n}\left(1+{\mathcal
O}\left(\frac{1}{\sqrt{n}}\right)\right)\\
&=\frac{\sqrt{2\sqrt{e}}}{4}n^{-1}\left(\frac{8}{e^2}\right)^n n^{2n}\left(1+{\mathcal
O}\left(\frac{1}{\sqrt{n}}\right)\right),
\end{align*}
where we use that $\sum_{k<\sqrt[4]{n}}\left(\frac{1}{k!2^{k}}\right)=
\sum^{\infty}_{k=0}\frac{1}{k!2^k} -\sum_{k\geq \sqrt[4]{n}}\left(\frac{1}{k!2^{k}}\right) =\sqrt{e} + \mathcal{O}(\frac{1}{\sqrt{n}})$.

For the second sum, we use the fact that $\binom{n}{n-k}N_{n+1}^{(n-k)}$ decreases as $k$  increases (Corollary~\ref{prop_III}) and again (\ref{asymp-binom-Nk}):
\begin{align*}
\Sigma_2&\leq n\cdot\binom{n}{\lceil\sqrt[4]{n}\rceil}N_{n+1}^{(n-\lceil\sqrt[4]{n}\rceil)}\\
&=\mathcal{O}\left(\frac{n}{\lceil\sqrt[4]{n}\rceil!2^{\sqrt[4]{n}}}\cdot n^{-1}\left(\frac{8}{e^2}\right)^n n^{2n}\right)\\
&=\mathcal{O}\left(n^{-3/2}\left(\frac{8}{e^2}\right)^n n^{2n}\right).
\end{align*}
Combining the above two estimates, we obtain that
\[
\mathrm{OGN}_n=\frac{\sqrt{2\sqrt{e}}}{4}n^{-1}\left(\frac{8}{e^2}\right)^n n^{2n}\left(1+{\mathcal
O}\left(\frac{1}{\sqrt{n}}\right)\right)
\]
which is the claimed result.\qed

\begin{Rem}\label{proof-poi-law}
Note that (\ref{asymp-binom-Nk}) (which gives the asymptotics of $\mathrm{OGN}_{n,n-k}$) and the result we just proved implies already the claim from Corollary~\ref{first-cor}-(i).
\end{Rem}

\section{Asymptotic Analysis of $\mathbf{GN_n}$}\label{enum-galled}

We now turn to the proof of Theorem~\ref{main-thm-1}. Our general strategy is to find upper and lower bounds for $\mathrm{GN}_n$ which admit the same first-order asymptotics. For this, we will use formula (\ref{form-GNn}) and some of the asymptotic tools from the last section.

We start with the following lemma.
\begin{lmm}
Define
\begin{equation}\label{Un}
U_n:=\sum_{{\mathcal T}}\prod_{v\in{\mathcal I}({\mathcal T})}\mathrm{OGN}_{c(v)},
\end{equation}
where the notation is as in Theorem~\ref{main-thm-GuRaZh}. Then, $\mathrm{GN}_n\leq U_n$.
\end{lmm}
\pf This follows immediately from Theorem~\ref{main-thm-GuRaZh} by noting that
\[
\binom{c_{\mathrm{lf}}(v)}{j-c_{\mathrm{nlf}}(v)}=\binom{c_{\mathrm{lf}}(v)}{c(v)-j}\leq\binom{c(v)}{c(v)-j}=\binom{c(v)}{j}
\]
and Proposition~\ref{gun-one-comp} and (\ref{rel-nkt}).\qed

Next, define the exponential generating function
\[
U(z):=\sum_{n\geq 0}U_n\frac{z^n}{n!}
\]
which can be obtained by recursive enumeration starting from the root of $\mathcal{T}$.

\begin{lmm}\label{lag-setting}
$U(z)$ satisfies the equations
\[
U(z)(1-M(U(z)))=z
\]
with
\[
M(z):=\sum_{\ell\geq 1}\frac{\mathrm{OGN}_{\ell+1}}{(\ell+1)!}z^{\ell}.
\]
Consequently,
\[
U_n=(n-1)![z^{n-1}](1-M(z))^{-n}.
\]
\end{lmm}
\pf In order to prove the first claim, observe that each tree $\mathcal{T}$ in (\ref{Un}) either consists of a single root, or a root to which two subtrees are attached, or a root to which three subtrees are attached, etc. Moreover, if $\ell$ subtrees are attached, then we have to give the root a weight of $\mathrm{OGN}_{\ell}$ to obtain (\ref{Un}). Overall, this gives
\[
U(z)=z+\sum_{\ell\geq 2}\frac{1}{\ell!}\cdot\mathrm{OGN}_{\ell}\cdot U(z)^{\ell},
\]
where in the sum, the term $1/\ell!$ is because the order of the subtrees is irrelevant, the term $\mathrm{OGN}_{\ell}$ is the weight and $U(z)^{\ell}$ is the exponential generating function of the $\ell$ subtrees. From this, the claimed equation for $U(z)$ follows.

Finally, the claimed expression for $U_n$ follows by Lagrange's inversion formula.\qed

To derive the first-order asymptotics of $U_n$, we use the following result of Bender and Richmond.

\begin{thm}[Bender and Richmond \cite{BeRi}]\label{BeRi}
Let $S(z)$ be a power series with $s_0=0, s_1\ne 0$ and $ns_{n-1}\sim\gamma s_n$. Then, for $\alpha\ne 0$ and $\beta$ real numbers, we have
\[
[z^n](1+S(z))^{\alpha n+\beta}\sim\alpha e^{\alpha a_1\gamma}ns_n.
\]
\end{thm}

Now, we can show the following.
\begin{pro}\label{main-thm-step-1}
We have, as $n\rightarrow\infty$,
\[
U_n\sim\frac{\sqrt{2e\sqrt[4]{e}}}{4}n^{-1}\left(\frac{8}{e^2}\right)^n n^{2n}.
\]
\end{pro}
\pf First, observe that from Proposition~\ref{main-pro-1} and Stirling's formula,
\[
[z^{\ell}]M(z)\sim\frac{2\sqrt[4]{e}}{\sqrt{\pi}}\ell^{-1/2}\left(\frac{8}{e}\right)^{\ell}{\ell}^{\ell}.
\]
From this we see that the sequences $[z^{\ell}]M(z)$ satisfies the assumptions from Theorem~\ref{BeRi} with $\gamma=1/8$. The claimed result follows now from that theorem and Stirling's formula.\qed

\begin{figure}[h]
\begin{center}
\begin{tikzpicture}
\draw (0cm,0cm) node[line width=0.03cm,inner sep=0,minimum size=0.2cm,draw,circle] (1) {};
\draw (-1.5cm,-1cm) node[line width=0.02cm,inner sep=0,minimum size=0.2cm,draw,circle] (2) {};
\draw (-0.5cm,-1cm) node[line width=0.02cm,inner sep=0,minimum size=0.2cm,draw,circle] (3) {};
\draw (0.5cm,-1cm) node[line width=0.02cm,inner sep=0,minimum size=0.2cm,draw,circle] (4) {};
\draw (1.5cm,-1cm) node[line width=0.02cm,inner sep=0,minimum size=0.2cm,draw,circle] (5) {};
\draw (-1.7cm,-2cm) node[line width=0.02cm,inner sep=0,minimum size=0.2cm,draw,circle] (6) {};
\draw (-1.3cm,-2cm) node[line width=0.02cm,inner sep=0,minimum size=0.2cm,draw,circle] (7) {};
\draw (-0.7cm,-2cm) node[line width=0.02cm,inner sep=0,minimum size=0.2cm,draw,circle] (8) {};
\draw (-0.3cm,-2cm) node[line width=0.02cm,inner sep=0,minimum size=0.2cm,draw,circle] (9) {};
\draw (-1cm,-1cm) node (10) {\small $\ldots$};
\draw (0.95cm,-1cm) node (11) {\small $\ldots$};

\draw[line width=0.02cm] (1)--(2);
\draw[line width=0.02cm] (1)--(3);
\draw[line width=0.02cm] (1)--(4);
\draw[line width=0.02cm] (1)--(5);
\draw[line width=0.02cm] (2)--(6);
\draw[line width=0.02cm] (2)--(7);
\draw[line width=0.02cm] (3)--(8);
\draw[line width=0.02cm] (3)--(9);
\draw [thick,decorate,decoration={brace,amplitude=4pt,mirror},xshift=0.4pt,yshift=-0.4pt](-1.9cm,-2.2cm)--(-0.1cm,-2.2cm) node[black,midway,yshift=-0.4cm] {\footnotesize $2j$};
\draw [thick,decorate,decoration={brace,amplitude=4pt,mirror},xshift=0.4pt,yshift=-0.4pt](0.3cm,-1.2cm)--(1.7cm,-1.2cm) node[black,midway,yshift=-0.4cm] {\footnotesize $n-2j$};
\end{tikzpicture}
\end{center}
\vspace*{-0.5cm}\caption{The phylogenetic trees $\mathcal{T}$ used in the definition of $L_n$ (labels of leaves are removed).}\label{star-cherry}
\end{figure}
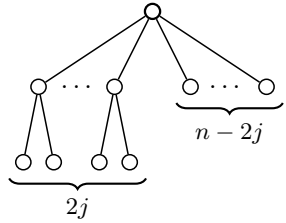

We next turn to the lower bound. Therefore, consider trees $\mathcal{T}$ in (\ref{form-GNn}) which consist of a root with children some of which are leaves and to the others we attach a cherry; see Figure~\ref{star-cherry}. We denote the number of galled networks with $n$ leaves arising from these trees $\mathcal{T}$ in (\ref{form-GNn}) by $L_n$. Clearly, $L_n\leq\mathrm{GN}_n$.

\vspace*{0.35cm}

For $L_n$, we have the following formula.

\begin{lmm}\label{form-Ln}
We have,
\begin{equation}\label{exp-Ln}
L_n=\sum_{j=0}^{\lfloor n/2\rfloor}\binom{n}{2j}\frac{(2j)!3^j}{j!}\sum_{\ell=0}^{n-2j}\binom{n-2j}{\ell}N_{n-j+1}^{(\ell+j)}.
\end{equation}
\end{lmm}
\pf Assume that exactly $j$ of the children of the root of $\mathcal{T}$ are followed by a cherry; see Figure~\ref{star-cherry}. Then, the outdegree of the root is $n-j$ and the contribution of these $\mathcal{T}$ to (\ref{form-GNn}) is
\[
\left(\sum_{\ell=0}^{n-2j}\binom{n-2j}{\ell}N_{n-j+1}^{(\ell+j)}\right)\left(\sum_{\ell=0}^{2}\binom{2}{\ell}N_{3}^{(\ell)}\right)^j.
\]
The second sum is $6$. Moreover, note that the number of $\mathcal{T}$ with the above property is
$\binom{n}{2j}\frac{(2j)!}{2^j j!}.$
Finally, summing over $j$ gives the claimed result.\qed
\\

Next, we consider the inner sum in (\ref{exp-Ln}).
\begin{lmm}\label{est-inner}
\begin{itemize}
\item[(i)] Uniformly for $0\leq j\leq\lfloor n/2\rfloor$, as $n\rightarrow\infty$,
\begin{equation}\label{est-1}
\sum_{\ell=0}^{n-2j}\binom{n-2j}{\ell}N_{n-j+1}^{(\ell+j)}={\mathcal O}\left(c^n n^{2n-2j}\right),
\end{equation}
where $c$ is a suitable constant.
\item[(ii)] We have, as $n\rightarrow\infty$,
\[
\sum_{\ell=0}^{n-2j}\binom{n-2j}{\ell}N_{n-j+1}^{(\ell+j)}=\frac{\sqrt{2\sqrt{e}}}{2^{3j+2}}n^{-1}\left(\frac{8}{e^2}\right)^nn^{2n-2j}\left(1+{\mathcal O}\left(\frac{j^2}{n}+\frac{1}{\sqrt{n}}\right)\right)
\]
uniformly for $j=o(\sqrt{n})$.
\end{itemize}
\end{lmm}
\pf We start by proving (i). First, by (\ref{ind-Nki}), we have
\[
N_{n-j+1}^{(\ell+j)}=\left(\prod_{s=0}^{\ell+j-1}\left(n-j-1+s+\frac{s}{2(n-j+s)}+{\mathcal O}\left(\frac{1}{n-j}\right)\right)\right)(2(n-j)-3)!!.
\]
Since $0\leq\ell\leq n-2j$, we have
\[
\prod_{s=0}^{\ell+j-1}\left(n-j-1+s+\frac{s}{2(n-j+s)}+{\mathcal O}\left(\frac{1}{n-j}\right)\right)={\mathcal O}(c_1^n n^{n-j})
\]
for a suitable constant $c_1$. Clearly, $(2(n-j)-3)!!={\mathcal O}(c_2^nn^{n-j})$ for a suitable constant $c$. Multiplying the last two estimates and plugging the product into (\ref{est-1}) gives
\[
\sum_{\ell=0}^{n-2j}\binom{n-2j}{\ell}N_{n-j+1}^{(\ell+j)}={\mathcal O}\left((2c_1c_2)^n n^{2n-2j}\right).
\]
Setting $c=2c_1c_2$ gives the claimed result.

We next show (ii). Here, from (\ref{asymp-binom-Nk}),
\[
\binom{n-j}{\ell}N_{n-j+1}^{(n-j-\ell)}=\frac{\sqrt{2}}{\ell!2^{\ell+3j+2}\sqrt[4]{e}}\left(\frac{8}{e^2}\right)^n n^{2n-2j}\left(1+{\mathcal O}\left(\frac{\ell^2+j^2+1}{n}\right)\right)
\]
uniformly for $\ell=o(\sqrt{n})$ and $j=o(\sqrt{n})$. Consequently, by Stirling's formula,
\begin{equation}\label{asymp-joint}
\binom{n-2j}{\ell}N_{n-j+1}^{(n-j-\ell)}=\frac{\sqrt{2}}{\ell!2^{\ell+3j+2}\sqrt[4]{e}}\left(\frac{8}{e^2}\right)^n n^{2n-2j}\left(1+{\mathcal O}\left(\frac{\ell^2+j^2+1}{n}\right)\right)
\end{equation}
uniformly for $\ell=o(\sqrt{n})$ and $j=o(\sqrt{n})$. From this, the claimed result follows with similar arguments as used at the end of the proof of Proposition~\ref{main-pro-1}.\qed

Now, we can derive the first-order asymptotics of $L_n$.
\begin{pro}\label{main-thm-step-2}
We have, as $n\rightarrow\infty$,
\[
L_n\sim\frac{\sqrt{2e\sqrt[4]{e}}}{4}n^{-1}\left(\frac{8}{e^2}\right)^nn^{2n}.
\]
\end{pro}

\pf We break the first sum in the formula for $L_n$ from Lemma~\ref{form-Ln} into two parts according to whether $j<\sqrt[4]{n}$ or not:
\begin{align*}
L_n&=\sum_{j<\sqrt[4]{n}}\binom{n}{2j}\frac{(2j)!3^j}{j!}\sum_{\ell=0}^{n-2j}\binom{n-2j}{\ell}N_{n-j+1}^{(\ell+j)}\\
&\qquad\quad+\sum_{\sqrt[4]{n}\leq j\leq\lfloor n/2\rfloor}\binom{n}{2j}\frac{(2j)!3^j}{j!}\sum_{\ell=0}^{n-2j}\binom{n-2j}{\ell}N_{n-j+1}^{(\ell+j)}=:L_n^{(1)}+L_{n}^{(2)}.
\end{align*}

For $L_n^{(2)}$, we use the uniform bound from Lemma~\ref{est-inner}-(i) and obtain that
\[
L_{n}^{(2)}={\mathcal O}\left(c^n\sum_{\sqrt[4]{n}\leq j\leq\lfloor n/2\rfloor}\frac{n!}{(n-2j)!}\cdot\frac{3^j}{j!}n^{2n-2j}\right)={\mathcal O}\left(nc^n\frac{3^{\lceil\sqrt[4]{n}\rceil}}{\lceil\sqrt[4]{n}\rceil!}n^{2n}\right).
\]

For $L_n^{(1)}$, using the expansion from Lemma~\ref{est-inner}-(ii) yields
\begin{align*}
L_n^{(1)}&=\left(\sum_{j<\sqrt[4]{n}}\frac{1}{j!}\left(\frac{3}{8}\right)^j\right)\frac{\sqrt{2\sqrt{e}}}{4}n^{-1}\left(\frac{8}{e^2}\right)^nn^{2n}
\left(1+\frac{1}{\sqrt{n}}\right)\\
&=\frac{\sqrt{2e\sqrt[4]{e}}}{4}n^{-1}\left(\frac{8}{e^2}\right)^nn^{2n}
\left(1+\frac{1}{\sqrt{n}}\right).
\end{align*}

Putting the two estimate together gives the claimed result.\qed

The proof of Theorem~\ref{main-thm-1} now directly follows from Proposition~\ref{main-thm-step-1} and Proposition~\ref{main-thm-step-2}.

\section{The Number of Reticulation Nodes}\label{num-of-ret}

In this section, we prove Theorem~\ref{main-thm-2} and the two corollaries from the introduction.

We first need a refinement of the sequence $L_n$ from the last section. Denote by $L_{n,k,j}$ the number of galled networks with $n$ leaves, $k$ reticulation nodes and $j$ inner reticulation nodes which arise again from trees whose root has children some of which are followed by cherries and some of which are leaves. Then, we have the following formula.

\begin{lmm}\label{form-Lnkj}
We have,
\[
L_{n,k,j}=\binom{n}{2j}\frac{(2j)!}{2^jj!}\sum_{\ell_1+\cdots+\ell_{j+1}=k-j}\left(\prod_{s=0}^{j}\binom{2}{\ell_s}N_3^{(\ell_s)}\right)\binom{n-2j}{\ell_{j+1}}N_{n-j+1}^{(j+\ell_{j+1})},
\]
where the sum over all non-negative integers with $\ell_s\in\{0,1,2\}$ for $0\leq s\leq j$.
\end{lmm}
\pf Note that galled networks counted by $L_{n,k,j}$ arise from the trees ${\mathcal T}$ depicted in Figure~\ref{star-cherry} whose number equals
\[
\binom{n}{2j}\frac{(2j)!}{2^jj!};
\]
compare with the proof of Lemma~\ref{form-Ln}.

What is left is to generate $k-j$ reticulation nodes that are followed by leaves; the sum in the claimed formula takes care of all the possibilities of picking these reticulation nodes from the $j$ cherries ($\ell_1,\ldots,\ell_j$) and the $n-2j$ leaves ($\ell_{j+1}$). Then, the internal nodes of ${\mathcal T}$ have to be replaced by the respective one-component galled networks and we are done.\qed

We can now prove Theorem~\ref{main-thm-2}.

\vspace*{0.35cm}\noindent{\it Proof of Theorem~\ref{main-thm-2}.} Recall that
\[
{\mathbb P}(X_n=j,n-Y_n=k)=\frac{\mathrm{GN}_{n,n-k,j}}{\mathrm{GN}_{n}};
\]
compare with (\ref{joint-dist}). Also, from the proof of Theorem~\ref{main-thm-1} in the last section, we know that, as $n\rightarrow\infty$,
\[
\mathrm{GN}_{n,n-k,j}=L_{n,n-k,j}+o(\mathrm{GN}_n)
\]
for all fixed $k$ and $j$. Thus, we need the asymptotics of $L_{n,n-k,j}$ as $n\rightarrow\infty$ for fixed $j$ and $k$.

In order to derive this asymptotics, first by Lemma~\ref{form-Lnkj}, we have
\[
L_{n,n-k,j}=\frac{1}{2^jj!}\cdot\frac{n!}{(n-2j)!}\sum_{\ell=\max\{0,j-k\}}^{2j}\left(\sum_{\ell_1+\cdots+\ell_j=\ell}\prod_{s=0}^{j}\binom{2}{\ell_s}N_3^{(\ell_s)}\right)\binom{n-2j}{\ell+k-j}N_{n-j+1}^{(n-k-\ell)}
\]
for $n\geq k+2j$. Next, note that
\[
\sum_{\ell_1+\cdots+\ell_j=\ell}\prod_{s=0}^{j}\binom{2}{\ell_s}N_3^{(\ell_s)}=[z^{\ell}]\left(\sum_{\ell=0}^{2}\binom{2}{\ell}N_3^{(\ell)}\right)^j=[z^{\ell}]\left(1+2z+3z^2\right)^j
\]
and by (\ref{asymp-joint}),
\[
\binom{n-2j}{\ell+k-j}N_{n-j+1}^{(n-k-\ell)}\sim\frac{\sqrt{2}}{(\ell+k-j)!2^{\ell+k+2j+2}\sqrt[4]{e}}\left(\frac{8}{e^{2}}\right)^n n^{2n-2j}.
\]
Consequently,
\[
L_{n,n-k,j}\sim\frac{\sqrt{2}}{2^{4j+2}j!\sqrt[4]{e}}\left(\sum_{\ell=\max\{0,j-k\}}^{2j}\frac{[z^{\ell}]\left(1+2z+3z^2\right)^j}{(\ell+k-j)!2^{\ell+k-j}}\right)\left(\frac{8}{e^2}\right)^n n^{2n}.
\]
Finally,
\[
\sum_{\ell=\max\{0,j-k\}}^{2j}\frac{[z^{\ell}]\left(1+2z+3z^2\right)^j}{(\ell+k-j)!2^{\ell+k-j}}=[z^{j-k}]e^{1/(2z)}(1+2z+3z^2)^j.
\]

Now, by putting everything together, we obtain that
\[
{\mathbb P}(X_n=j,n-Y_n=k)\sim\frac{L_{n,n-k,j}}{\mathrm{GN}_n}\sim\frac{e^{-7/8}}{16^j j!}[z^{j-k}]e^{1/(2z)}(1+2z+3z^2)^j,
\]
where we used Theorem~\ref{main-thm-1}. This proves the claimed result.\qed

What is left is to prove the corollaries.

\vspace*{0.35cm}\noindent{\it Proof of Corollary~\ref{first-cor}.} Part (i) follows from (\ref{asymp-binom-Nk}); compare with Remark~\ref{proof-poi-law}. Alternatively, we can use Theorem~\ref{main-thm-2} since
\begin{align*}
{\mathbb P}(n-Z_n=k)={\mathbb P}(n-Y_n=k\vert X_n=0)&=\frac{{\mathbb P}(X_n=0, n-Y_n=k)}{{\mathbb P}(X_n=0)}\\
&\longrightarrow\frac{{\mathbb P}(X=0, Y=k)}{{\mathbb P}(X=0)}.
\end{align*}
By a simple computation
\[
\frac{{\mathbb P}(X=0, Y=k)}{{\mathbb P}(X=0)}=\frac{e^{-1/2}}{2^kk!}
\]
which proves the claimed result.

As for part (ii), observe that the limit law of $X_n$ is given by
\[
{\mathbb P}(X=j)=\frac{e^{-7/8}}{16^jj!}\sum_{k\geq-j}[z^{j-k}]e^{1/(2z)}(1+2z+3z^2)^j.
\]
Now,
\[
\sum_{k\geq-j}[z^{j-k}]e^{1/(2z)}(1+2z+3z^2)^j=\sum_{k\leq 2j}[z^k]e^{1/(2z)}(1+2z+3z^2)^j=e^{1/2}6^j.
\]
Consequently,
\[
{\mathbb P}(X=j)=\frac{e^{-3/8}}{j!}\left(\frac{3}{8}\right)^j
\]
which proves the claim from part (ii).\qed

\vspace*{0.35cm}\noindent{\it Proof of Corollary~\ref{second-cor}.} Theorem~\ref{main-thm-2} implies that
\[
{\mathbb E}(n-Y_n)\longrightarrow{\mathbb E}(Y)\qquad\text{and}\qquad{\rm Var}(n-Y_n)={\rm Var}(Y_n)\longrightarrow{\rm Var}(Y).
\]
So, what we have to do is to evaluate the mean and variance of $Y$.

For the mean, we have
\begin{equation}\label{mean-Y}
{\mathbb E}(Y)=\sum_{j\geq 0}\frac{e^{-7/8}}{16^jj!}\sum_{k\geq -j}k[z^{j-k}]e^{1/(2z)}(1+2z+3z^2)^j.
\end{equation}
The second sum can be rewritten as follows
\[
\sum_{k\geq -j}k[z^{j-k}]e^{1/(2z)}(1+2z+3z^2)^j=\sum_{k\leq 2j}(j-k)[z^{k}]e^{1/(2z)}(1+2z+3z^2)^j.
\]
Now, note that
\[
\sum_{k\leq 2j}j[z^{k}]e^{1/(2z)}(1+2z+3z^2)^j=e^{1/2}j6^j
\]
and
\[
\sum_{k\leq 2j}k[z^{k}]e^{1/(2z)}(1+2z+3z^2)^j=\frac{{\rm d}}{{\rm d} z}e^{1/(2z)}(1+2z+3z^2)^j\Big\vert_{z=1}=e^{1/2}\left(\frac{4}{3}j-\frac{1}{2}\right)6^j.
\]
Overall,
\[
\sum_{k\geq -j}k[z^{j-k}]e^{1/(2z)}(1+2z+3z^2)^j=e^{1/2}\left(\frac{1}{2}-\frac{j}{3}\right)6^j.
\]
Plugging this into (\ref{mean-Y}) and straightforward simplification gives
\[
{\mathbb E}(Y)=\frac{3}{8}.
\]

For the variance of $Y$, a similar computation proves the second claim.\qed

\section{Asymptotically Counting Dup-Trees}\label{dup-trees}

As was pointed out in \cite{GuRaZh}, one-component galled networks are in close relationship with leaf-multi-labeled trees (or LML trees for short). In this section, we will recall this relationship and present results which either directly follow from our results for one-component galled networks or are obtained with a similar method of proof.

We start by recalling some definitions. First, a (binary, rooted) {\it LML tree} is a leaf-labeled tree with labels of the leaves not necessarily distinct. An LML tree is called {\it dup-tree} if each label can be used at most twice. Obviously,  binary phylogenetic trees are dup-trees, where label repetition is prohibited.

A {\it cherry} of a tree is a pair of leaves that are adjacent to a common non-leaf node. If the two leaves have the same label, we call the cherry a {\it twin-cherry}. A dup-tree is called {\it twin-cherry free} if it does not contain a twin-cherry.

\begin{figure}[h]
\begin{center}
\begin{tikzpicture}
\draw (0cm,0cm) node[inner sep=1.2pt,line width=0.8pt,draw,circle] (1) {{\footnotesize $\rho$}};
\draw (0cm,-0.8cm) node[minimum size=7.8pt,line width=0.8pt,draw,circle] (2) {};
\draw (-0.5cm,-1.6cm) node[minimum size=7.8pt,line width=0.8pt,draw,circle] (3) {};
\draw (1cm,-3cm) node[minimum size=7.8pt,line width=0.8pt,draw,circle,fill=gray!35] (6) {};
\draw (-1cm,-2.4cm) node[minimum size=7.8pt,line width=0.8pt,draw,circle] (4) {};
\draw (0cm,-2.4cm) node[minimum size=7.8pt,line width=0.8pt,draw,circle] (5) {};
\draw (-1.5cm,-3.2cm) node[inner sep=1.2pt,line width=0.8pt,draw,circle] (7) {{\footnotesize $2$}};
\draw (-0.5cm,-3.2cm) node[minimum size=7.8pt,line width=0.8pt,draw,circle,fill=gray!35] (8) {};
\draw (-0.5cm,-4cm) node[inner sep=1.2pt,line width=0.8pt,draw,circle] (9) {{\footnotesize $1$}};
\draw (1cm,-3.8cm) node[inner sep=1.2pt,line width=0.8pt,draw,circle] (10) {{\footnotesize $3$}};

\draw[-stealth,line width=0.8pt] (1) -- (2);
\draw[-stealth,line width=0.8pt] (2) -- (3);
\draw[-stealth,line width=0.8pt] (2) -- (6);
\draw[-stealth,line width=0.8pt] (3) -- (4);
\draw[-stealth,line width=0.8pt] (3) -- (5);
\draw[-stealth,line width=0.8pt] (4) -- (8);
\draw[-stealth,line width=0.8pt] (5) -- (8);
\draw[-stealth,line width=0.8pt] (5) -- (6);
\draw[-stealth,line width=0.8pt] (4) -- (7);
\draw[-stealth,line width=0.8pt] (8) -- (9);
\draw[-stealth,line width=0.8pt] (6) -- (10);

\draw[stealth-stealth,line width=0.8pt] (2.5cm,-2cm) -- (4cm,-2cm);

\draw (6.5cm,-0.8cm) node[inner sep=2.5pt,line width=0.8pt,draw,circle] (2) {};
\draw (6cm,-1.6cm) node[inner sep=2.5pt,line width=0.8pt,draw,circle] (3) {};
\draw (7.5cm,-3.2cm) node[inner sep=0.8pt,line width=0.8pt,draw,circle] (6) {{\tiny $3$}};
\draw (5.5cm,-2.4cm) node[inner sep=2.5pt,line width=0.8pt,draw,circle] (4) {};
\draw (6.5cm,-2.4cm) node[inner sep=2.5pt,line width=0.8pt,draw,circle] (5) {};
\draw (5.2cm,-3.2cm) node[inner sep=0.8pt,line width=0.8pt,draw,circle] (7) {{\tiny $2$}};
\draw (5.8cm,-3.2cm) node[inner sep=0.8pt,line width=0.8pt,draw,circle] (8) {{\tiny $1$}};
\draw (6.2cm,-3.2cm) node[inner sep=0.8pt,line width=0.8pt,draw,circle] (9) {{\tiny $1$}};
\draw (6.8cm,-3.2cm) node[inner sep=0.8pt,line width=0.8pt,draw,circle] (10) {{\tiny $3$}};

\draw[line width=0.8pt] (2) -- (3);
\draw[line width=0.8pt] (2) -- (6);
\draw[line width=0.8pt] (3) -- (4);
\draw[line width=0.8pt] (3) -- (5);
\draw[line width=0.8pt] (4) -- (8);
\draw[line width=0.8pt] (5) -- (9);
\draw[line width=0.8pt] (5) -- (10);
\draw[line width=0.8pt] (4) -- (7);
\end{tikzpicture}
\end{center}
\caption{The bijection between one-component galled networks and twin-cherry free dup-trees.}\label{1-GN-to-dup}
\end{figure}
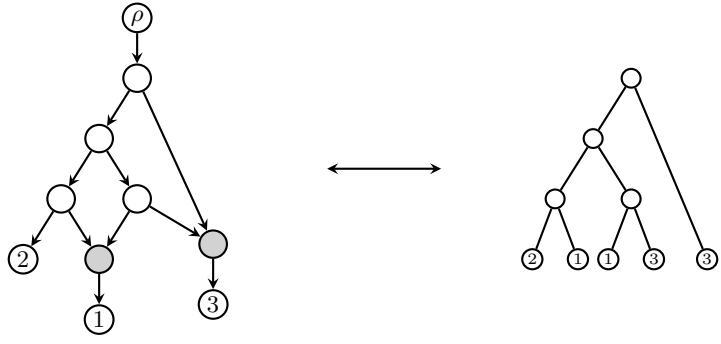

\begin{pro}[Gunawan et al. \cite{GuRaZh}]
There is a bijection between one-component galled networks with $n$ leaves and $k$ reticulation nodes and twin-cherry free dup-trees with $n$ different labels exactly $k$ of which are repeated.
\end{pro}
The bijection is actually easy to construct: remove the pendant edge below a reticulation node and replace the reticulation node by two labeled leaves with the label of the removed leaf. Then, attach these two leaves to the parents of the removed reticulation node; see Figure~\ref{1-GN-to-dup}.

\vspace*{0.35cm}

Denote by $\mathrm{FDU}_n$ the number of twin-cherry free dup-trees with $n$ distinct labels. Then, by Proposition~\ref{main-pro-1} and Corollary~\ref{first-cor}-(i), we have the following result.
\begin{thm}\label{dup-tree-1}
For the number $\mathrm{FDU}_n$ of twin-cherry free dup-trees with $n$ distinct labels, as $n\rightarrow\infty$,
\[
\mathrm{FDU}_n\sim\frac{\sqrt{2\sqrt{e}}}{4}n^{-1}\left(\frac{8}{e^2}\right)^nn^{2n}.
\]
Moreover, the number $\mathrm{FR}_n$ of repeated labels of a twin-cherry free dup-tree with $n$ distinct labels picked uniformly at random satisfies the following limit distribution result:
\[
n-\mathrm{FR}_n\stackrel{d}{\longrightarrow}\mathrm{Poi}(1/2),\qquad (n\rightarrow\infty).
\]
\end{thm}

In fact, we can find the first-order asymptotics of the number $\mathrm{DU}_n$ of all (not necessarily twin-cheery free) dup-trees with $n$ distinct labels as well. Therefore, we recall the following recursive way for computing this number, which was stated in the conclusion of \cite{GuRaZh} (compare with Proposition~\ref{gun-one-comp}).

\begin{pro}[Gunawan et al. \cite{GuRaZh}]
The number $\mathrm{DU}_n$ of dup-trees with $n$ distinct leaves is given by
\[
\mathrm{DU}_n=\sum_{k=0}^{n}\binom{n}{k}B_{n+1}^{(k)},
\]
where $B_{n}^{(k)}$ is recursively given by
\[
B_{n}^{(k)}=(n+k-2)B_{n}^{(k-1)}+\frac{1}{2}\sum_{1\leq d\leq k-1}\binom{k-1}{d}(2d-1)!!\left(B_{n-d}^{(k-1-d)}-B_{n-d+1}^{(k-1-d)}\right)
\]
for $2\leq k\leq n-1$ with initial values $B_{n}^{(0)}=(2n-5)!!$ and $B_n^{(1)}=(n-1)(2n-5)!!$.
\end{pro}

From this proposition, with the same method of proof as in Section~\ref{enum-1-c}, we obtain the following result.
\begin{thm}\label{enum-dup}
For the number $\mathrm{DU}_n$ of dup-trees with $n$ distinct labels, as $n\rightarrow\infty$,
\[
\mathrm{DU}_n\sim\frac{\sqrt{2e\sqrt{e}}}{4}n^{-1}\left(\frac{8}{e^2}\right)^nn^{2n}.
\]
Moreover, the number $\mathrm{R}_n$ of repeated labels of a dup-tree with $n$ distinct labels picked uniformly at random satisfies the following limit distribution result:
\[
n-\mathrm{R}_n\stackrel{d}{\longrightarrow}\mathrm{Poi}(1/2),\qquad (n\rightarrow\infty).
\]
\end{thm}

In fact, instead of re-doing the analysis from Section~\ref{enum-1-c}, one can alternatively use what we have already proved for $N_{n}^{(k)}$ by exploring the following relation between $\mathrm{OGN}_{n,k}$ and $\mathrm{DU}_n$ as well as its refinement for $\mathrm{OGN}_{n,k}$ and $\mathrm{DU}_{n,k}$, where the latter denotes the number of dup-trees with $n$ distinct labels exactly $k$ of which are repeated.

\begin{lmm}
We have,
\[
\mathrm{DU}_{n}=\sum_{k=0}^{n}2^{n-k}\times\mathrm{OGN}_{n,k}
\]
and
\[
\mathrm{DU}_{n,k}=\sum_{\ell=0}^{k}\binom{n-\ell}{k-\ell}\times\mathrm{OGN}_{n,\ell}
\]
\end{lmm}

\pf Recall that $\mathrm{OGN}_{n,k}$ also counts the number of twin-cherry free dup-trees with $n$ distinct labels exactly $k$ of which are repeated. Now, the claimed results follow by observing that leaves with labels which are not repeated can be either replaced by a twin-cherry or left unchanged.\qed

We now use this to prove Theorem~\ref{enum-dup}.

\vspace*{0.3cm}\noindent{\it Proof of Theorem~\ref{enum-dup}.} From (\ref{asymp-binom-Nk}), we obtain that
\[
2^{k}\times\mathrm{OGN}_{n,n-k}=\frac{\sqrt{2}}{k!4\sqrt[4]{e}}n^{-1}\left(\frac{8}{e^2}\right)^n n^{2n}\left(1+{\mathcal O}\left(\frac{k^2+1}{n}\right)\right)
\]
uniformly in $k=o(\sqrt{n})$. Then, with similar arguments as in the last paragraph of the proof of Proposition~\ref{main-pro-1},
\begin{align*}
\mathrm{DU}_{n}=\sum_{k=0}^{n}2^k\times\mathrm{OGN}_{n,n-k}&\sim\frac{\sqrt{2}}{4\sqrt[4]{e}}\left(\sum_{k\geq 0}\frac{1}{k!}\right)n^{-1}\left(\frac{8}{e^2}\right)^n n^{2n}\\
&=\frac{\sqrt{2e\sqrt{e}}}{4}n^{-1}\left(\frac{8}{e^2}\right)^n n^{2n}
\end{align*}
which is the claimed result for $\mathrm{DU}_n$.

Next, for the distribution of $R_n$, observe that
\[
\mathrm{DU}_{n,n-k}=\sum_{\ell=0}^{k}\binom{k+\ell}{\ell}\times\mathrm{OGN}_{n,n-k-\ell}.
\]
Again, from (\ref{asymp-binom-Nk}),
\[
\binom{k+\ell}{\ell}\times\mathrm{OGN}_{n,n-k-\ell}=\frac{\sqrt{2}}{\ell!k!2^{k+l+2}\sqrt[4]{e}}n^{-1}\left(\frac{8}{e^2}\right)^n n^{2n}\left(1+{\mathcal O}\left(\frac{\ell^2+1}{n}\right)\right)
\]
uniformly for $\ell=o(\sqrt{n})$. Thus,
\[
\mathrm{DU}_{n,n-k}\sim\frac{1}{k!2^k}\cdot\frac{\sqrt{2\sqrt{e}}}{4}n^{-1}\left(\frac{8}{e^2}\right)^n n^{2n}.
\]
which implies the claimed result for the distribution of $R_n$ since
\[
{\mathbb P}(R_n=k)=\frac{\mathrm{DU}_{n,n-k}}{\mathrm{DU}_n}.
\]
This concludes the proof of Theorem~\ref{enum-dup}.\qed

Finally, Theorem~\ref{dup-tree-1} and Theorem~\ref{enum-dup} imply the following corollary.
\begin{cor}
The fraction of twin-cherry-free dup-trees with $n$ distinct labels amongst all dup-trees with $n$ distinct labels tends to $e^{-1/2}$ as $n$ tends to infinity.
\end{cor}

\section{Conclusion}\label{con}

In this paper, we have derived the first-order asymptotics of the number of galled networks and proved limit laws for shape parameters of galled networks which are picked uniformly at random. This is the first time that such results are obtained for a widely-used class of phylogenetic networks of ``large" size; compare with the discussion in Section~\ref{intro}.

We end the paper with some concluding remarks.

First,  because of (\ref{rel-nkt}), Theorem~\ref{main-thm-2} also implies corresponding limit distribution results for the number of tree nodes and for the total number of nodes of galled networks with $n$ leaves. It would be interesting to study stochastic properties of the height of galled networks and other shape parameters, where the height is defined as the number of edges on each of the longest paths from the root to a leaf (see the conclusion of \cite{DiSeWe} where this parameter was called the {\it depths}).

Second,  another question is to investigate further properties of the limit distribution $Y$ from Theorem~\ref{main-thm-2} since the expression given there (which is obtained by summing over $j$) is not particularly easy to handle. (However, as seen in Section~\ref{num-of-ret}, at least the computation of moments is feasible from this expression.)

Finally, the most immediate question raised by our study is how about other ``large" classes of phylogenetic networks? Can they be (asymptotically) enumerated as well? Also, how to study the limit behavior of shape parameters for them? One natural class to consider next would be the class of reticulation-visible networks. In \cite{GuRaZh}, the authors asked for (exact) enumeration results. We now broaden this question and ask for an asymptotic study similar to the one we carried out in this paper.

\section*{Acknowledgement}

We thank Hsien-Kuei Hwang for drawing our attention to \cite{BeRi} which was the key for proving our asymptotic result for $\mathrm{GN}_n$. We also thank the (anonymous) reviewer for a careful reading.

\newpage

\section*{Appendix}

\vspace{0.5cm}
\begin{table}[bp!]
\begin{center}
\rotatebox{90}{\begin{minipage}{21cm}\centering
{\small\begin{tabular}{c|cccccccccc}
$k\setminus n$& 2 & 3 & 4 & 5 & 6 & 7 & 8 & 9 & 10 & 11 \\
\hline
0 & 1 & 1 & 3 & 15 & 105 & 945 & 10,395 & 135,135 & 2,027,025 & 34,459,425 \\
1 & 0 & 1 & 6 & 45 & 420 & 4,725 & 62,370 & 945,945 & 16,216,200 & 310,134,825 \\
2 & - & 3 & 20 & 189 & 2,160 & 28,875 & 442,260 & 7,640,325 & 147,026,880 & 3,119,591,475 \\
3 & - & - & 87 & 993 & 13,407 & 207,135 & 3,603,915 & 69,757,065 & 1,487,243,835 & 34,639,019,415 \\
4 & - & - & - & 6,249 & 97,182 & 1,701,855 & 33,121,890 & 709,428,825 & 16,587,636,030 & 420,498,508,815 \\
5 & - & - & - & - & 804,585 & 15,738,765 & 338,588,685 & 7,946,584,695 & 202,099,078,125 & 5,537,451,658,725 \\
6 & - & - & - & - & - & 161,685,045 & 3,808,469,970 & 97,162,333,695 & 2,669,506,204,050 & 78,595,220,899,125 \\
7 & - & - & - & - & - & - & 46,726,507,485 & 1,287,228,175,056 & 37,987,475,258,565 & 1,195,779,444,849,675 \\
8 & - & - & - & - & - & - & - & 18,363,976,595,055 & 579,247,192,040,580 & 19,410,597,807,225,345 \\
9 & - & - & - & - & - & - & - & - & 9,420,991,174,195,965 & 334,803,875,697,765,495 \\
10 & - & - & - & - & - & - & - & - & - & 6,114,381,201,716,874,975
\end{tabular}}
\caption{The values of $N_{n}^{(k)}$ for $2\leq n\leq 11$ and $0\leq k\leq n$.}\label{n-n-k}

\vspace{2cm}
\begin{tabular}{c|cccccccc}
$j\setminus k$ & 0 & 1 & 2 & 3 & 4 & 5 & 6 & 7 \\
\hline
0 & 46,726,507,485 & 26,659,289,790 & 7,110,362,385 & 1,159,266,150 & 125,137,025 & 9,287,460 & 436,590 & 10,395 \\
1 & 18,868,231,935 & 20,820,564,765 & 12,078,633,735 & 3,747,731,400 & 692,176,275 & 79,858,170 & 5,554,395 & 186,795 \\
2 & 4,976,625,150 & 7,604,859,780 & 5,995,908,765 & 2,779,284,375 & 813,268,575 & 145,143,495 & 14,794,920 & 686,700 \\
3 & 960,639,750 & 1,795,456,530 & 1,708,006,230 & 983,507,175 & 366,209,550 & 86,543,100 & 11,981,970 & 746,235 \\
4 & 122,089,275 & 260,763,300 & 281,838,690 & 186,377,625 & 80,515,575 & 22,424,850 & 3,717,000 & 281,925 \\
5 & 7,577,955 & 17,681,895 & 20,896,785 & 15,181,425 & 7,243,425 & 2,242,485 & 416,115 & 35,595
\end{tabular}
\caption{The values of $\mathrm{GN}_{7,7+j-k,j}$ with $0\leq j\leq 5$ and $0\leq k\leq 7$.}\label{ll-ret}
\end{minipage}}
\end{center}
\end{table}
\end{document}